\documentclass[a4paper,10pt]{amsart}

\usepackage{epsfig}
\usepackage{amsthm,amsfonts}
\usepackage{amssymb,graphicx,color}
\usepackage{float}
\graphicspath{ {Images/} }
\usepackage[labelfont={bf},textfont={it}]{caption}
\usepackage[labelfont={bf, it}]{subcaption}
\usepackage[all]{xy}
\usepackage{verbatim}
\usepackage{hyperref}
\usepackage{enumitem}

\newtheorem{theorem}{Theorem}[section]
\newtheorem*{theorem*}{Theorem}
\newtheorem{lemma}[theorem]{Lemma}
\newtheorem{corollary}[theorem]{Corollary}
\newtheorem{proposition}[theorem]{Proposition}
\newtheorem{definition}[theorem]{Definition}

\newtheorem{claim}{Claim}
\newtheorem{example}[theorem]{Example}
\newtheorem{remark}[theorem]{Remark}

\newcommand{\R}{\mathbb{R}}

\newcommand{\Z}{\mathbb{Z}}
\newcommand{\Q}{\mathbb{Q}}
\newcommand{\N}{\mathbb{N}}

\begin{document}

\title[Global bi-Lipschitz classification of semialgebraic surfaces]
{Global bi-Lipschitz classification of semialgebraic surfaces}

\author[A. Fernandes]{Alexandre Fernandes}
\author[J. E. Sampaio]{Jos\'e Edson Sampaio}

\address{Alexandre Fernandes: Departamento de Matem\'atica, Universidade Federal do Cear\'a, Av. Humberto Monte, s/n Campus do Pici - Bloco 914, 60455-760, Fortaleza-CE, Brazil. E-mail: {\tt alex@mat.ufc.br}}
\address{J. Edson Sampaio: Departamento de Matem\'atica, Universidade Federal do Cear\'a, Rua Campus do Pici, s/n, Bloco 914, Pici, 60440-900, Fortaleza-CE, Brazil. E-mail: {\tt edsonsampaio@mat.ufc.br}}

\thanks{The first named author was partially supported by CNPq-Brazil grant 304700/2021-5. The second author was partially supported by CNPq-Brazil grant 310438/2021-7.
This study was financed in part by the CAPES-BRASIL Finance Code 001.}

\keywords{bi-Lipschitz classification, surfaces, Nash surfaces}
\subjclass[2010]{58A07; 14R05; 14P25; 14B05; 32S50}

\begin{abstract}
We classify semialgebraic surfaces in $\mathbb{R}^n$ with isolated singularities up to bi-Lipschitz homeomorphisms with respect to the inner distance. In particular, we obtain complete classifications for the Nash surfaces and the complex algebraic curves. We also address the minimal surfaces with finite total curvature.
\end{abstract}

\maketitle

\section{Introduction}
In 1999, in the seminal paper on bi-Lipschitz classification of 2D real singularities \cite{Birbrair}, Lev Birbrair proved the following result

\begin{theorem}[Theorem of Birbrair]\label{thm:birbrair}
	Given the germ of a semialgebraic set, $(X,a)$, with isolated singularity and connected link, there is a unique rational number $\beta\geq 1$ such that $(X,a)$ is bi-Lipschitz homeomorphic, with respect its inner distance, to the  germ at $0\in\R^3$ of the $\beta$-horn
	$$\{(x,y,z)\in\R^3 \ \colon \ x^2+y^2=z^{2\beta} \ \mbox{and} \ z\geq 0\}.$$ 
\end{theorem}
A similar result was also obtained in \cite{Grieser:2003}.

The goal of this present paper is to bring the Theorem of Birbrair and its ideas to a global perspective on the inner bi-Lipschitz geometry of the 2D real subsets in $\R^n$. Such a perspective, in the smooth case, is closely related to Fu Conjecture (see \cite{Fu:1998}), which states that {\it a complete Riemannian surface in $\R^3$ with $\int K^+<2\pi$ and $\int K^-<+\infty$ and which is homeomorphic to $\R^2$ must be bi-Lipschitz homeomorphic to $\R^2$}. 
This conjecture was positively answered by Bonk and Lang in \cite{BonkL:2003} for surfaces endowed with the inner distance. Another related study was presented in \cite{BelenkiiB:2005} by Belen'ki\u{\i} and Burago, they presented a classification of complete Aleksandrov surfaces with finite total curvature under some restrictions in their singularities and such that their ends have non-zero growth speed.

Coming back to our goal, let us start by recalling the topological classification of compact (without boundary) smooth surfaces. It is well-known (since the 1860s) that, given a compact smooth surface $S$ in $\R^n$, two symbols $\theta_S\in \{-1,1\}$ and $g_S\in\N\cup \{0\}$ complete determine $S$ up to diffeomorphisms, where $\theta_S$ says that $S$ is orientable or not; and $g_S$ is the genus of $S$. In the setting of (not necessarily compact) properly embedded smooth surfaces in $\R^n$, in order to have some control on the topology of such surfaces, let us assume they are semialgebraic. In some sense, since compact manifolds (without boundary) are diffeomorphic to semialgebraic ones (see \cite{Nash}), that assumption is not too restrictive. In this setting, there is a topological  structure theorem that says: $\exists$ a radius $R>0$ such that for any $\rho \geq R$, the Euclidean sphere $\mathbb{S}_{\rho}$ intersects transversally $S$ and  $S\setminus B_{\rho}$ the set of points in $S$ and outside the Euclidean ball $B_{0,\rho}$ is diffeomorphic to the cylinder $[\rho,\infty)\times S\cap\mathbb{S}_{\rho}$. The connected components of $S\setminus B_{\rho}$ are called the {\it ends} of $S$, and any two family of ends corresponding to $\rho$ and $\rho'$, respectively, are always diffeomorphic each other, for any pairs of radius $\rho,\rho'>R$ (see \cite{Coste}). 

Then, in the case of properly embedded smooth surfaces $S$ in $\R^n$ which are semialgebraic, so-called {\it Nash surfaces in $\R^n$}, we have a list of three symbols to determine $S$ up to  diffeomorphism, namely: $\theta_S\in \{-1,1\}$, $g_S\in\N\cup \{0\}$ and $e_S\in\N\cup \{0\}$ that is the number of ends of $S$. 

In this paper, we consider Nash surfaces in $\R^n$ equipped with the inner distance
$$d_{inn}(x_1,x_2) = \inf\{length(\gamma) \ \colon \ \gamma \ \mbox{is a path on} \ S \ \mbox{connecting} \ x_1,x_2\in S \}$$
and we classify those surfaces up to bi-Lipschitz homeomorphisms with respect to the inner distance, the so-called {\it inner lipeomorphims}. Actually, associated to each Nash surface $S$, we present a list of symbols, $\theta_S\in \{-1,1\}$, $g_S\in\N\cup \{0\}$, $e_S\in\N\cup \{0\}$ and $\beta_1$, ...,$\beta_{e_S}$, where $\beta_{i}'$s ($\leq 1$) are rational numbers associated to the ends of $S$; which determines $S$ up to inner lipeomorphisms.

Finally, we address semialgebraic surfaces with isolated inner Lipschitz singularities, we classify all these surfaces up inner lipeomorphims by using a combinatorial invariant so called inner Lipschitz code (see Definition \ref{def:inner_code}), we also bring some applications of this classification to complex algebraic plane curves and minimal surfaces with finite total curvature.

\section{Preliminaries}

Given a path connected subset $X\subset\R^n$, the
\emph{inner distance}  on $X$  is defined as follows: given two points $x_1,x_2\in X$, $d_{X,inn}(x_1,x_2)$  is the infimum of the lengths of paths on $X$ connecting $x_1$ to $x_2$. 

\begin{definition}\label{def:lne}
Let $X\subset\R^n$ be a subset. We say that $X$ is {\bf Lipschitz normally embedded (LNE)} if there exists a constant $c\geq 1$ such that $d_{X,inn}(x_1,x_2)\leq C\|x_1-x_2 \|$, for all pair of points $x_1,x_2\in X$. 
\end{definition}

For instance, considering the real (resp. complex) cusp $x^2=y^3$, in $\R^2$ (resp. in $\mathbb{C}^2$), one can see that this set is not LNE.

\begin{definition}\label{lipschitz function}
Let $X\subset\R^n$ and $Y\subset\R^m$. A mapping $f\colon X\rightarrow Y$ is called {\bf outer} (resp. {\bf inner}) {\bf Lipschitz} if there exists $\lambda >0$ such that is
$$\|f(x_1)-f(x_2)\|\le \lambda \|x_1-x_2\| \quad (\mbox{resp. } d_{X,inn}(f(x_1),f(x_2))\le \lambda d_{X,inn}(x_1,x_2))$$  for all $x_1,x_2\in X$. A outer Lipschitz (resp. inner Lipschitz) mapping $f\colon X\rightarrow Y$ is called
{\bf outer} (resp. {\bf inner}) {\bf lipeomorphism} if its inverse mapping exists and is outer Lipschitz (resp. inner Lipschitz) and, in this case, we say that  $X$ and $Y$ are {\bf outer} (resp. {\bf inner}) {\bf lipeomorphic}.
\end{definition}

Let $X\subset\R^n$ be a closed 2-dimensional semialgebraic set in the following two definitions,   

\begin{definition}
	A point $p\in X$ is called {\bf topologically regular} if there exists a neighborhood $V\subset X$ of $p$ homeomorphic to an open disc in $\R^2$. When all the points $p\in X$ are topologically regular, $X$ is said a {\bf semialgebraic topological surface} in $\R^n$. 
\end{definition}

\begin{definition}
	A point $p\in X$ is called {\bf inner Lipschitz regular} if there exists a neighborhood $V\subset X$ of $p$ inner lipeomorphic to an open disc in $\R^2$; otherwise it is called  {\bf inner Lipschitz singular}. We denote by ${\rm Reg}_{inLip}(X)$ (resp. ${\rm Sing}_{inLip}(X)$) the set of all inner Lipschitz regular (resp. singular) points of $X$.
\end{definition}

\subsection{Semialgebraic ends, topological and Lipschitz singularities}

Let $S$ be a closed semialgebraic subset of $\R^n$. Let us assume that $S$ has possible isolated singularities, i.e., there exists a finite subset $\Sigma\subset S$ such that all points in $S\setminus\Sigma$ are inner Lipschitz regular points of $S$. Such a subset is what we call a {\bf semialgebraic surface in $\R^n$ with isolated inner Lipschitz singularities}. As a consequence of the Local Conic Structure Theorem (by using the inversion mapping $z\mapsto z/|z|^2$), we see that there exists a large radius $R>0$ such that 
\begin{enumerate}
	\item $S$ is transversal to the Euclidean sphere $\mathbb{S}(0,\rho)$ for any $\rho\geq R$.
	\item there exists a semialgebraic homeomorphism $$\phi\colon S\setminus B(0,R)\rightarrow \{ t\cdot u \ : \ t\geq R \ \mbox{and} \ u\in \mathbb{S}(0,R)\}$$ such that $|\phi (z)|= |z|$ for any $z\in S$ outside of the Euclidean ball $B(0,R)$.
\end{enumerate}
It follows that $S\setminus B(0,R)$ has finite many semialgebraic connected components $S_1,\dots S_{e_S}$, each $S_i$ is semialgebraicly homeomorphic to the cylinder $\mathbb{S}^1\times[R,\infty)$. Moreover, each $S_i$ is  semialgebraicly homeomorphic to $S_i\setminus B(0,\rho)$ for any $\rho\geq R$. Those subsets $S_1,\dots S_{e_S}$ are called the {\bf ends} of $S$. Notice that the ends of $S$ are well-defined up to semialgebraic homeomorphisms. When $S$ has only one end, we say that it is {\bf connected at infinity}.

Next, we recall the notion of tangent cone at infinity which is important to the study of Lipschitz geometry of ends of semialgebraic sets.

\begin{definition}
Let $X\subset \R^m$ be an unbounded subset. We say that $v\in \R^m$ is {\bf tangent to $X$ at infinity} if there are a sequence of real positive numbers $\{ t_j \}_{j\in \N}$ such that $t_j\to +\infty$ and a sequence of points $\{x_j\}_{j\in \N}\subset X$ such that $\lim\limits _{j\to +\infty }\frac{1}{t_j}x_j=v$. Denote by $C(X, \infty )$ the set of $v\in\R^m$ which are tangent to  $X$ at infinity and we call it the {\bf tangent cone of $X$ at infinity}.
\end{definition}

\subsection{Contact of curves at infinity}

\begin{definition}
	Let $\Gamma_1, \Gamma_2\subset \R^n$ be two unbounded semialgebraic curves, which are connected at infinity. Fixed $K>1$, we define $f_{\Gamma_1,\Gamma_2}^K\colon (0,+\infty)\to \R$ by  
	$$
	f_{\Gamma_1,\Gamma_2}^K(r)=dist(A^K_r(\Gamma_1),A^K_r(\Gamma_2)),
	$$
	where $A^K_r(X)=\{y\in X;\frac{r}{K}\leq \|y\|\leq Kr\}$ and $dist(X,Y)=\inf\{\|x-y\|;x\in X$ and $y\in Y\}$. If $\Gamma_1\cap  \Gamma_2$ is an unbounded set, we define $Cont^K(\Gamma_1, \Gamma_2)=-\infty$ and if $\Gamma_1\cap  \Gamma_2$ is a bounded set, we define 
	$$
	Cont^K(\Gamma_1, \Gamma_2)=\lim\limits_{r\to +\infty}\frac{\log{f_{\Gamma_1,\Gamma_2}^K(r)}}{\log{r}}.
	$$
\end{definition}

\begin{remark}\label{prop:contact_vs_tangency}
	Let $\Gamma_1, \Gamma_2\subset \R^n$ be unbounded semialgebraic curves, which are connected at infinity. Let $K>1$. Then $Cont^K(\Gamma_1, \Gamma_2)=1$ if and only if $C(\Gamma_1,\infty)\not=C(\Gamma_2,\infty)$.
\end{remark}

\begin{proposition}\label{prop:non-dependency_of_K}
	Let $\Gamma_1, \Gamma_2\subset \R^n$ be unbounded semialgebraic curves, which are connected at infinity. Let $K,\tilde K>1$. Then $Cont^K(\Gamma_1, \Gamma_2)=Cont^{\tilde K}(\Gamma_1, \Gamma_2)$.
\end{proposition}
\begin{proof}
	It follows from the definition that $Cont^K(\Gamma_1, \Gamma_2)=-\infty$ if and only if $Cont^{\tilde K}(\Gamma_1, \Gamma_2)=-\infty$. So, we may assume that $Cont^K(\Gamma_1, \Gamma_2)$ and $Cont^{\tilde K}(\Gamma_1, \Gamma_2)$ are finite numbers.
	From Remark \ref{prop:contact_vs_tangency}, we may also assume that $C(\Gamma_1,\infty)=C(\Gamma_2,\infty)$.

	We assume that $K<\tilde K$. 
	
	\begin{claim}\label{claim:Ksquare}
		If $\tilde K\leq K^2$ then $Cont^K(\Gamma_1, \Gamma_2)=Cont^{\tilde K}(\Gamma_1, \Gamma_2)$.
	\end{claim}
	\begin{proof}[Proof of Claim \ref{claim:Ksquare}]
		Since $K<\tilde K$, it is clear that $f_{\Gamma_1,\Gamma_2}^{\tilde K}(r)\leq f_{\Gamma_1,\Gamma_2}^K(r)$. Since $C(\Gamma_1,\infty)=C(\Gamma_2,\infty)$, we obtain 
		$$
		\min\{f_{\Gamma_1,\Gamma_2}^K(Kr),f_{\Gamma_1,\Gamma_2}^K(r/K)\}\lesssim f_{\Gamma_1,\Gamma_2}^{\tilde K}(r) \mbox{ as }r\to +\infty.
		$$
		Moreover, 
		$$
		f_{\Gamma_1,\Gamma_2}^K(r)\approx f_{\Gamma_1,\Gamma_2}^K(r/K) \mbox{ and }f_{\Gamma_1,\Gamma_2}^K(r)\approx f_{\Gamma_1,\Gamma_2}^K(r/K) \mbox{ as }r\to +\infty,
		$$ 
		and thus we obtain $f_{\Gamma_1,\Gamma_2}^K(r)\approx f_{\Gamma_1,\Gamma_2}^{\tilde K}(r)$ as $r\to +\infty$, which gives 
		$$
		Cont^K(\Gamma_1, \Gamma_2)=Cont^{\tilde K}(\Gamma_1, \Gamma_2).
		$$
		
	\end{proof}
	It follows from Claim \ref{claim:Ksquare} that $Cont^K(\Gamma_1, \Gamma_2)=Cont^{K^m}(\Gamma_1, \Gamma_2)$ for all positive integer $m$. Let $m$ be an positive integer such that $K^m\leq \tilde K\leq K^{2m}$. By Claim \ref{claim:Ksquare} again, $Cont^{K^m}(\Gamma_1, \Gamma_2)=Cont^{\tilde K}(\Gamma_1, \Gamma_2)$, which finishes the proof.
\end{proof}

Thus, we define $Cont(\Gamma_1, \Gamma_2)=Cont^K(\Gamma_1, \Gamma_2)$ for some $K>1$.

\begin{proposition}\label{prop:contact_invariance}
	Let $\Gamma_1, \Gamma_2\subset \R^n$ and  $\tilde\Gamma_1, \tilde\Gamma_2\subset \R^m$ be unbounded semialgebraic curves, which are connected at infinity. Assume that there exists an outer lipeomorphism $F\colon \Gamma_1\cup \Gamma_2\to \tilde\Gamma_1\cup \tilde\Gamma_2$ such that $F(\Gamma_i)=\tilde \Gamma_i$, $i=1,2$. Then $Cont(\Gamma_1, \Gamma_2)=Cont(\tilde\Gamma_1, \tilde\Gamma_2)$.
\end{proposition}
\begin{proof}
	Since $F$ is an outer lipeomorphism, there is $M\geq 1$ such that  
	$$
	\frac{1}{M}\|x-y\|\leq \|F(x)-F(y)\|\leq M\|x-y\|, \quad \forall x,y\in  \Gamma_1\cup \Gamma_2.
	$$
	Let $x_0\in  \Gamma_1\cup \Gamma_2$ such that $\|x_0\|\geq 1$. Then, for any $x\in  \Gamma_1\cup \Gamma_2$ such that $\|x\|\geq r_0=\max\{3\|x_0\|,3M\|F(x_0)\|\}$, we have
	\begin{eqnarray*}
		\|F(x)\|&\leq & \|F(x)-F(x_0)\|+\|F(x_0)\| \\
		&\leq &M\|x-x_0\|+\|x\|\\
		&\leq &3M\|x\|
	\end{eqnarray*}
	and 
	\begin{eqnarray*}
		\|F(x)\|&\geq & \|F(x)-F(x_0)\|-\|F(x_0)\| \\
		&\geq & \frac{1}{M}\|x-x_0\|-\frac{1}{3M}\|x\|\\
		&\geq & \frac{1}{M}\|x\|-\frac{1}{M}\|x_0\|-\frac{1}{3M}\|x\|\\
		&\geq & \frac{1}{3M}\|x\|.
	\end{eqnarray*}
	
	Therefore, for any $K>1$ and $\tilde K= 3MK$, we have that $F(A^K_r(\Gamma_1\cup \Gamma_2))\subset A^{\tilde K}_r(\tilde \Gamma_1\cup \tilde \Gamma_2)$ for all $r\geq Kr_0$.
	Thus, 
	$$
	dist(A^K_r(\Gamma_1),A^K_r(\Gamma_2))\geq \frac{1}{M}dist(F(A^K_r(\Gamma_1)),F(A^K_r(\Gamma_2))) \geq dist(A^{\tilde K}_r(\tilde \Gamma_1),A^{\tilde K}_r(\tilde \Gamma_2))
	$$
	for all $r\geq Kr_0$, which implies $Cont(\Gamma_1, \Gamma_2)\geq Cont(\tilde\Gamma_1, \tilde\Gamma_2)$.
	Similarly, we also prove that $Cont(\Gamma_1, \Gamma_2)\leq Cont(\tilde\Gamma_1, \tilde\Gamma_2)$, which finishes the proof.
\end{proof}

\begin{example}\label{exam:contact_horn}
	Let $\beta\in \mathbb{Q}$ with $\beta\leq 1$. Let $\Gamma_1=\{(x,0)\in\R^2;x\geq 1\}$ and $\Gamma_2=\{(x,y)\in \R^2;x\geq 1$ and $y=x^{\beta}\}$. Then $Cont(\Gamma_1, \Gamma_2)=\beta$.
\end{example}

\subsection{Lipeomorphisms between circles}\label{sec:circle_lipeomorphisms}

The main goal of this Subsection is to show that, if $f,g\colon\mathbb{S}^1\rightarrow\mathbb{S}^1$ are two lipeomorphisms with the same orientation, then there exists a lipeotopy $H_t\colon\mathbb{S}^1\rightarrow\mathbb{S}^1$ ($0\leq t\leq 1$) such that $H_0=f$ and $H_1=g$. By lipeotopy we mean a lipeomorphism $H\colon [0,1]\times\mathbb{S}^1\rightarrow[0,1]\times\mathbb{S}^1$ of the type $H(t,x)=(t,H_t(x))$, which is equivalent to the following: $H_t$ is as a family of lipeomorphisms with uniform constant. Possibly this result is already known, but we did not find an appropriate reference to quote, this is why we present a proof of it.

Let $P\colon\R\rightarrow\mathbb{S}^1$ be the covering mapping $P(x)=e^{2\pi ix}$.

\begin{lemma}
	If $\phi\colon\mathbb{S}^1\rightarrow\mathbb{S}^1$ is a positive homeomorphism such that $\phi(1)=1$, then there exists a unique positive homeomorphism $\widetilde{\phi}\colon\R\rightarrow\R$ such that $\widetilde{\phi}(0)=0$ and $P\circ\widetilde{\phi}=\phi\circ P$ (in particular, $\widetilde{\phi}(x+ n)=\widetilde{\phi}(x)+n$ $\forall n\in\Z$). Conversely, for each positive homeomorphism $\widetilde{\phi}\colon\R\rightarrow\R$ such that $\widetilde{\phi}(0)=0$ and $\widetilde{\phi}(x+ n)=\widetilde{\phi}(x)+n$ $\forall n\in\Z$, there is a unique positive homeomorphism  $\phi\colon\mathbb{S}^1\rightarrow\mathbb{S}^1$ such that $\phi(1)=1$ and $P\circ\widetilde{\phi}=\phi\circ P$. Finally, $\phi$ is a lipeomorphism iff $\widetilde{\phi}$ is a lipeomorphism. 
\end{lemma} 

\begin{proof} 
	Let $\phi\colon\mathbb{S}^1\rightarrow\mathbb{S}^1$ is a positive homeomorphism such that $\phi(1)=1$. Let $\widetilde{\phi}\colon\R\rightarrow\R$ be defined by: given $x\in\R$, let $\gamma\colon [0,x]\rightarrow\mathbb{S}^1$ be the path defined by $\gamma(t)=\phi(e^{2\pi i t})$ and  let $\widetilde{\gamma}\colon [0,x]\rightarrow\R$ be the lifting of $\gamma$ by the covering mapping $P$ with $\widetilde{\gamma}(0)=0$; so, $\widetilde{\phi}(x):=\widetilde{\gamma}(x)$. By definition, we have $\widetilde{\phi}(0)=0$ and $P\circ\widetilde{\phi}=\phi\circ P$, and, since $\phi$ is positive and $\widetilde{\phi}$ is a local homeomorphism, $\widetilde{\phi}$ is an increasing homeomorphism from $\R$ to $\R$. 
	
	Conversely, let $\widetilde{\phi}\colon\R\rightarrow\R$ be an increasing  homeomorphism such that $\widetilde{\phi}(0)=0$ and $\widetilde{\phi}(x+ n)=\widetilde{\phi}(x)+n$ $\forall n\in\Z$. Then, $\phi\colon\mathbb{S}^1\rightarrow\mathbb{S}^1$ defined by $\phi(e^{2\pi ix})=e^{2\pi\widetilde{\phi}(x)}$ is a positive homeomorphism such that $\phi(1)=1$.
	
	Finally, we are going to show that $\phi$ is a lipeomorphism iff $\widetilde{\phi}$ is a lipeomorphism with the same constants. Let us consider $\R$ and $\mathbb{S}^1$ equipped with the standard Riemannian Metric. Thus, $P\colon\R\rightarrow\mathbb{S}^1$ comes as a local (on each interval  of length $2\pi$) isometry, hence $\widetilde{\phi}$ lipeomorphism implies that $\phi$ is a lipeomorphism with the same constants. From the other hand, if $\phi$ is a lipeomorphism with constants $c\geq 1$ and $1/c$, we have that $\widetilde{\phi}$ is a lipeomorphism with these constants on each interval of length $2\pi$. Now, given $a<b$ in $\R$, we partition the interval $[a,b]$ with subintervals of length smaller than $2\pi$: $a=x_0<x_1<\cdots<x_{n-1}<x_n=b$, and:
	
	\medskip
	
	\begin{eqnarray*}
	|\widetilde{\phi}(b)-\widetilde{\phi}(a)| &=& \widetilde{\phi}(b)-\widetilde{\phi}(a) \\
	&=& \sum_{j=1}^{n} \widetilde{\phi}(x_j)-\widetilde{\phi}(x_{j-1})\\
	&\leq& \sum_{j=1}^{n} c(x_j-x_{j-1}) \\
	&=& c|b-a|
	\end{eqnarray*}
	and
	\begin{eqnarray*}
	|\widetilde{\phi}(b)-\widetilde{\phi}(a)| &=& \widetilde{\phi}(b)-\widetilde{\phi}(a) \\
	&=& \sum_{j=1}^{n} \widetilde{\phi}(x_j)-\widetilde{\phi}(x_{j-1})\\
	&\geq& \sum_{j=1}^{n} \frac{1}{c}(x_j-x_{j-1}) \\
	&=& \frac{1}{c}|b-a|
	\end{eqnarray*}
\end{proof}

Once we have the above lemma, given a positive lipeomorphim $\phi\colon\mathbb{S}^1\rightarrow\mathbb{S}^1$ such that $\phi(1)=1$, let us consider $\widetilde{H}_t\colon\R\rightarrow\R$ defined by $\widetilde{H}_t(x) = (1-t)\widetilde{\phi}(x)+tx.$ We see that $\widetilde{H}_t(0)=0$ and $\widetilde{H}_t$ is a family of positive lipeomorphisms (with uniform constant) such that.
 $$\widetilde{H}_t(x+n)=\widetilde{H}_t(x)+n \ \forall n\in\Z.$$ Then, $H_t$ given by the above  lemma is a family of lipeomorphisms from $\mathbb{S}^1$ to $\mathbb{S}^1$ (with uniform constant) such that $H_0=\phi$ and $H_1= id_{\mathbb{S}^1}$ 

\begin{proposition}\label{prop:lipeotopy}
	Let $f,g\colon\mathbb{S}^1\rightarrow\mathbb{S}^1$ be two lipeomorphisms with the same orientation. Then there exists a lipeotopy $H_t\colon\mathbb{S}^1\rightarrow\mathbb{S}^1$ such that $H_0=f$ and $H_1=g$.
\end{proposition}

\begin{proof}
	We do the proof in the case $f(1)=g(1)$, i.e. $f\circ g^{-1}(1)=1$. Then, by the  previous discussion,  there exists a lipeotopy $H_t\colon\mathbb{S}^1\rightarrow\mathbb{S}^1$ such that $H_0=f\circ g^{-1}$ and $H_1=id_{\mathbb{S}^1}$. Finally, we get $K_t:= H_t \circ g$ give us a lipeotopy such that $K_0=f$ and $K_1=g$.
\end{proof}

\section{Ends of semialgebraic surfaces in $\R^n$}

\subsection{Infinity strips}\label{sec:triangles}

Let $a>0$ and $\beta\in\Q$; $\beta\leq 1$. Let us denote $$T_{\beta}=\{(x,y)\in\R^2 \ \colon \ a\leq x \ \mbox{and} \ 0\leq y\leq x^{\beta}\} .$$ Notice that, up to outer lipeomorphims, the definition of $T_{\beta}$ does not depend on the constant $a>0$.

\begin{lemma}\label{lemma: lne}
	$T_{\beta}$ is LNE.
\end{lemma}

\begin{proof}
	Since $T_{\beta}$ is a convex subset of $\R^2$ in the case $\beta > 0$, we are going to prove this lemma for $\beta < 0$. Let $P,Q\in T_{\beta}$. The the length of segment $\overline{PQ}$ is exactly $|Q-P|$. So, if the segment $\overline{PQ}$ is contained in $T_{\beta}$, we have $d_{inn}(P,Q)=|Q-P|$, otherwise $\overline{PQ}$ intersects the boundary $\{(x,x^{\beta}) \ \colon \ 1\leq x \}$ of the set $T_{\beta}$ into two points $A=(a,a^{\beta})$ and $B=(b,b^{\beta})$ ($a<b$). 
\begin{claim}
	The length $l(\gamma)$ of the boundary path $\gamma\colon[a,b]\rightarrow T_{\beta}$; $\gamma(t)=(t,t^{\beta})$ is bounded by $2|B-A|$.
\end{claim}
In fact, once  $l(\gamma)=\int_{a}^{b}|\gamma'(t)|dt$, we have
\begin{eqnarray*}
	l(\gamma) &=& \int_{a}^{b}\sqrt{1+\beta^2t^{2\beta-2}}dt \\
	&\leq & \int_{a}^{b} (1-\beta t^{\beta-1})dt \quad \mbox{(see that $\beta < 0$)} \\
	&=& (b-a) - (b^{\beta}-a^{\beta}) \\
	&=& (b-a) + (a^{\beta}-b^{\beta}) \\
	&\leq& 2 |B-A|
\end{eqnarray*}

The claim is proved. 

\medskip

Finally, once we have proved the claim, we get
\begin{eqnarray*}
d_{inn}(P,Q) &\leq& |A-P|+l(\gamma)+|Q-B| \\
&\leq& |A-P|+2 |B-A|+|Q-B| \\
&\leq& 2 |Q-P|.
\end{eqnarray*}
This finishes the proof that $T_{\beta}$ is LNE.
\end{proof}

\begin{definition}\label{def:triangle}
	Let $X\subset\R^n$ be a semialgebraic subset. We say that $X$ is a {\bf $\beta$-strip at infinity} if there exist a compact subset $K\subset\R^n$ and a germ of a semialgebraic inner lipeomorphism $F\colon X\setminus K\rightarrow T_{\beta}$.
\end{definition}

\begin{remark}
	As an immediate consequence of Example \ref{exam:contact_horn}, we have no ambiguity to define $\beta$-strip at infinity, in other words, if $T_{\beta}$ is inner lipeomorphic to $T_{\beta'}$ then $\beta=\beta'$.
\end{remark}

\begin{proposition}\label{prop: triangle-f}
	Let $f\colon[a,\infty)\rightarrow\R$ be a positive semialgebraic function such that $f(x)\approx x^{\beta}$ as $x\to\infty$ for some rational number $\beta\leq 1$. In this case $$X=\{(x,y)\in\R^2 \ \colon \ a\leq x \ \mbox{and} \ 0\leq y\leq f(x)\}$$ is a LNE $\beta$-strip at infinity.
\end{proposition}

\begin{proof}
	By assumption, we have a real number $c>0$ and $f(x)=cx^{\beta}+o_{\infty}(x^{\beta})$ where, $\displaystyle\frac{o_{\infty}(x^{\beta})}{x^{\beta}}\to 0$ as $x\to\infty$. 
	
	\medskip
	
	Let $F\colon T_{\beta}\rightarrow X$ be defined by $\displaystyle 
	F(x,y)=(x,\frac{yf(x)}{cx^{\beta}})$. It is clear that $F$ is a semialgebraic homeomorphism. The Jacobian matrix $DF(x,y)$ is bounded as we see below
	$$
	DF(x,y)=\left(\begin{array}{ccc}
	1& &0\\
	& & \\
	\frac{y}{cx^{2\beta}}[f'(x)x^{\beta}-f(x)x^{\beta-1}]& & \frac{f(x)}{cx^{\beta}}\\
	\end{array}\right)
	$$
	and, also, its determinant is bounded and away from zero as $x\to \infty$. This proves that $F$ is an outer lipeomorphim which give us that $X$ is LNE and a $\beta$-strip at infinity.
\end{proof}

 Consider the following semialgebraic arcs on the boundary of $T_{\beta}$ 
 $$\gamma_1=\{(x,y)\in T_{\beta} \ \colon \ y=x^{\beta} \} \quad \mbox{and} \quad \gamma_2=\{(x,y)\in T_{\beta} \ \colon \ y=0 \}.$$ In the case of a $\beta$-strip at infinity $X$, its {\it boundary arcs at infinity} are $F(\gamma_1)$ and $F(\gamma_2)$ where  $F \colon X\setminus K\rightarrow T_{\beta}$ is any semialgebraic inner lipeomorphism and $K$ is a compact subset of $X$. 

\begin{lemma}\label{lemma:glue_2_triangles}
	Let $X_i$ be a $\beta_i$-strip at infinity, $i=1,2$. If $X_1\cap X_2$ is a common boundary arc to $X_1$ and $X_2$, then $X_1\cup X_2$ is a $\beta$-strip at infinity with $\beta=\max\{\beta_1,\beta_2\}$. 
\end{lemma}

\begin{proof}
	Let $F_1\colon X_1\rightarrow T_{\beta_1}$ and $F_2\colon X_2\rightarrow T_{\beta_2}^*$ be semialgebraic inner lipeomorphisms; where $T_{\beta}^*=\{(x,y)\in\R^2 \ \colon \  (x,-y)\in T_{\beta}\}$. For each $x\geq 1$, let us denote $r_i(x)=|F_i^{-1}(x,0)|$. We see that $r_i$ is a semialgebraic outer lipeomorphism function; and $R_i(x,y)=(r_i(x),y)$ give us a semialgebraic outer lipeomorphism ($i=1,2$). Then, we define
	$$
	F(z)=
	\begin{cases}
	R_1(F_1(z)) ,& \mbox{if} \ z\in X_1 \\
	R_2(F_2(z)) ,& \mbox{if} \ z\in X_2
	\end{cases}
	$$	
	
	Since $R_1(F_1(z))=R_2(F_2(z))$ for all $z\in X_1\cap X_2$, we have $F$ is a continuous and semialgebraic mapping. Now, we are going to show that $F$ is an inner Lipschitz mapping.
	
	In fact, we know that $F_|X_1$ and $F_|X_2$ are inner lipeomorphism, then there exists a constant $c$ such that
	$d_{inn}(F(z_1),F(z_2))\leq c d_{inn}(z_1,z_2)$ if $z_1,z_2\in X_i$ ($i=1,2$). Thus, let us consider the case $z_1\in X_1$ and $z_2\in X_2$. Let $\gamma$ be a path on $X_1\cup X_2$ connecting $z_1$ to $z_2$ such that $d_{inn}(z_1,z_2)=l(\gamma)$. Then, we can write $\gamma=\gamma_1 * \cdots *\gamma_r$ in such a way that each $\gamma_j$ is a path on $X_1$ or $X_2$. Let us denote by $a_j$ the initial point and $b_j$ the final point of $\gamma_j$. Thus,
	\begin{eqnarray*}
		d_{inn}(F(z_1),F(z_2)) &\leq& \sum d_{inn}(F(b_j),F(a_j)) \\
		&\leq& \sum cd_{inn}(b_j,a_j) \\
		&\leq& c\sum l(\gamma_j) \\
		&=& cl(\gamma) \\
		&=& cd_{inn}(z_1,z_2) 
	\end{eqnarray*}

This proves that $F$ is an inner Lipschitz mapping. So, we can use similar arguments to show that $F^{-1}$ is also an inner Lipschitz mapping.

\medskip

Finally, we are going to show that the image of $F$ is a $\beta$-strip at infinity. Since $x\mapsto r_i(x)$ is an outer lipeomorphim ($i=1,2$), we have the image $F(X_1\cup X_2)=R_1(T_{\beta})\cup R_2(T_{\beta}^*)$ is the following subset of $\R^2$
$$I=\{(x,y)\in\R^2 \ \colon \ x\geq 1 \ \mbox{and} \ -f_2(x)\leq y\leq f_1(x)\}$$ where $f_1,f_2\colon[1,\infty)\rightarrow\R $ are semialgebraic positive functions such that
$$ f_i(x) \approx x^{\beta_i} \quad \mbox{as} \quad x\to\infty \quad i=1,2.$$ Then, we see the mapping $(x,y)\mapsto (x,y+f_2(x))$ gives a semialgebraic outer lipeomorphism between the image $I$ and the set below
$$ J=\{(x,y)\in\R^2 \ \colon \ x\geq 1 \ \mbox{and} \ 0\leq y\leq f_1(x)+f_2(x)\}.$$ Thus, since $[f_1(x)+f_2(x)]\approx x^{\beta}$ as $x\to \infty$ ($\beta=\max\{\beta_1,\beta_2\}$), by Proposition \ref{prop: triangle-f}, it follows that  $J$ is a $\beta$-strip at infinity.

\end{proof}

\begin{proposition}[Gluing of Strips]\label{lemma:gluing_strips}
	Let $X_1,\dots,X_r$ be semialgebraic subsets of $\R^n$ such that:
	\begin{enumerate}
		\item[a)] $X_i$ is a $\beta_i$-strip at infinity, $i=1,\dots,r$.
		\item[b)] $X_i\cap X_{i+1}$ is a common boundary arc to $X_i$ and $X_{i+1}$, $i=1,\dots,r-1$.
		\item[c)] $X_i\cap X_j=\emptyset$ if $|i-j|>1$.
	\end{enumerate}
In this case, $X_1\cup\cdots\cup X_r$ is a $\beta$-strip at infinity, where $\beta=\max\{\beta_1,\dots,\beta_r\}$.
\end{proposition}

\begin{proof}
	It is imediate consequence from Lemma \ref{lemma:glue_2_triangles}.
\end{proof}
 
As an immediate consequence of the proof of Lemma \ref{lemma:glue_2_triangles}, we can state the following lemma.
\begin{lemma}[Parametrization Lemma]\label{lemma:parametrization}
	Let $X$ be a $\beta$-strip at infinity with boundary arcs $\gamma_1$ and $\gamma_2$. Then, there exist a compact subset $K\subset X$ and a semialgebraic inner lipeomorphism $F\colon X\setminus K\rightarrow T_{\beta}$ such that;
	
-  For $z\in\gamma_1\setminus K$, $F(z)=(|z|,|z|^{\beta})$, .

- For $z\in\gamma_2\setminus K$, $F(z)=(|z|,0)$.
\end{lemma}


\subsection{Tubes}\label{sec:tubes}

Given a rational number $\beta \leq 1$, let us denote
$$P_{\beta}=\{(x,y,z)\in\R^3 \ \colon \ x^2+y^2=z^{2\beta} \ \mbox{and} \ z\geq a\}$$ where $a>0$. It is important to mention that, up to outer lipeomorphims, the definition of $P_{\beta}$ does not depend on $a$.

\begin{definition}\label{def:tube}
	Let $X\subset\R^n$ be a semialgebraic subset. We say that $X$ is a {\bf $\beta$-tube} if there exist a compact subset $K\subset\R^n$ and a germ of a semialgebraic inner lipeomorphism $F\colon X\setminus K\rightarrow P_{\beta}$.
\end{definition}

\begin{remark}
	We have no ambiguity to define $\beta$-tube, in other words, if $P_{\beta}$ is inner lipeomorphic to $P_{\beta'}$ then $\beta=\beta'$.
\end{remark}

\begin{proposition}\label{prop:glue-tube}
	Let $X_1,\dots X_r$ be semialgebraic subsets of $\R^n$ such that:
	\begin{enumerate}
		\item[a)] $X_i$ is a $\beta_i$-strip at infinity, $i=1,\dots,r$.
		\item[b)] if $r=2$, then $X_1\cap X_2$ is the union of the boundary arcs of $X_1$ and $X_2$;
		\item[c)] if $r>2$, then $X_i\cap X_{i+1}$ is a common boundary arc to $X_i$ and $X_{i+1}$, $i=1,\dots,r$ (here, $X_{r+1}:=X_1$) and $X_i\cap X_j=\emptyset$ if $1<|i-j|<r$.
	\end{enumerate}
	In this case, $X_1\cup\cdots\cup X_r$ is a $\beta$-tube, where $\beta=\max\{\beta_1,\dots,\beta_r\}$.
	
\end{proposition}

\begin{proof}
	Without loss of generality, one may assume $\beta_1 = \beta$. Let us write $X_1$ as a union of two other $\beta$-strips at infinity $X_{1,1}$ and $X_{1,2}$ such that $X_{1,1}\cap X_{1,2}$ is a common boundary arc to $X_{1,1}$ and $X_{1,2}$. On may suppose that $X_{1,2}$ shares a boundary arc with $X_2$ and $X_{1,1}$ shares a boundary arc with $X_r$. So, the family $X_{1,2},X_2,\dots,X_r$ satisfies the conditions of Gluing of Strips in Proposition \ref{lemma:gluing_strips}, hence $Y_2=X_{1,2}\cup X_2\cdots\cap X_r$ is a $\beta$-strip at infinity. Then, $X$ is the union of two $\beta$-strip at infinity $Y_1=X_{1,1}$ and $Y_2$ ($Y_2$ defined above), such that $Y_1\cap Y_2$ is the union of the boundary arcs of $Y_1$ and $Y_2$. 
	
	Let us consider the following decomposition of $P_{\beta}$:
	$$P_{\beta}^1=\{(x,y,z)\in P_{\beta} \ \colon \ x\geq 0\} \ \mbox{and} \ P_{\beta}^2=\{(x,y,z)\in P_{\beta} \ \colon \ x\leq 0\}.$$ We see that $P_{\beta}^1$ and $P_{\beta}^2$ are $\beta$-strip at infinity, $P_{\beta}=P_{\beta}^1\cup P_{\beta}^2$, and $P_{\beta}^1\cap P_{\beta}^2$ is the union of the boundary arcs of $P_{\beta}$. Since $P_{\beta}^1$ and $P_{\beta}^2$ are $\beta$-strip at infinity, we have a compact subset of $\R^n$ and semialgebraic inner lipeomorphism $$ F_1\colon Y_1\setminus K\rightarrow P_{\beta}^1 \quad \mbox{and} \quad F_2\colon Y_2\setminus K\rightarrow P_{\beta}^2 $$ such that $|F(z)|=|z|$ for any $z$ belonging to the boundary arcs of $Y_i$, $i=1,2$ (see Parametrization Lemma \ref{lemma:parametrization}).
	
	Finally, the mapping $F\colon X\setminus K\rightarrow P_{\beta}$ defined by
	$$
	F(z)=
	\begin{cases}
	F_1(z) ,& \mbox{if} \ z\in Y_1\setminus K \\
	F_2(z) ,& \mbox{if} \ z\in Y_2\setminus K
	\end{cases}
	$$	
	is a semialgebraic inner lipeomorphism.
\end{proof}

\begin{theorem}\label{thm:beta_ends}
	Let $S\subset\R^n$ be a semialgebraic surface with isolated inner Lipschitz singularities. For each end of $S$, let us say $S_i$, there is a unique rational $\beta_i\leq 1$ such that $S_i$ is a $\beta_i$-tube.
\end{theorem}

In order to proof this theorem, we need to recall the notion of $L$-regular sets. Such subsets of $\R^n$ are defined by induction on $n$ (see \cite{KP}).  Given $x\in\R^n$, let us write $x=(x',x_n)\in\R^{n-1}\times\R$. A semialgebraic subset $X\subset\R^n$ is called a {\it standard L-regular cell} in $\R^n$, with constant $C>0$ if: $X=\{0\}$ for $n=0$, and for $n>0$ the set $X$ is of one of the following types:

\medskip

\noindent{(\it graph)} $$X=\{(x',x_n)\in\R^{n-1}\times\R \ \colon \ x_n=h(x'); \ x'\in X'\}$$

\medskip

\noindent{(\it band)}  $$X=\{(x',x_n)\in\R^{n-1}\times\R \ \colon \ f(x')<x_n<g(x'); \ x'\in X'\}$$

\noindent where $X'\subset\R^{n-1}$ is a standard L-regular cell in $\R^{n-1}$ with constant $C$, $f,g,h\colon X'\rightarrow\R$ are $C^1$ semialgebraic functions such that 
$$f(x')<g(x') \quad \forall \ x'\in X'$$ and $$ |df(x')|\leq C, \ |dg(x')|\leq C\ \mbox{and} \ \ |dh(x')|\leq C, \quad \forall \ x'\in X'.$$ 

In general, a semialgebraic subset $Z\subset\R^n$ is called a {\it L-regular cell} in $\R^n$, with constant $C>0$, if there exists an orthogonal change of variables $\Psi\colon\R^n\rightarrow\R^n$ such that $\Psi(Z)$ is a standard L-regular cell in $\R^n$ with constant $C$.
 
\begin{proposition}\label{prop:regular-cell}
	Let $X\subset\R^n$ be a 2-dimensional L-regular cell in $\R^n$ (with constant $C>0$). If $X$ is unbounded and has only one end, then $\overline{X}$ is a $\beta$-strip at infinity for some rational number $\beta\leq 1$.
\end{proposition}

\begin{proof}
	It is enough to assume $X$ is a standard L-regular cell in $\R^n$. This proof is by induction on $n$.  Since $X$ is 2-dimensional, we have $n\geq 2$. 

\medskip

\noindent{\it Case $n=2$}. In this case, necessarily $X$ is a band, let us say
$$X=\{(x_1,x_2)\in\R\times\R \ \colon \ f(x_1)<x_2<g(x_1), \ x_1\in X'\}$$ where $X'$ is an open interval in $\R$. Since $X$ is unbounded and $|df(x_1)|\leq C$ and $|dg(x_1)|\leq C$ for all $x_1\in X'$, we get $X'$ is also unbounded and has only one end. Let us suppose $X=(a,\infty)$. Thus, the closure $\overline{X}$ is outer lipeomorphic to the set
$$\{(x_1,x_2)\in\R\times\R \ \colon \ 0\leq x_2 \leq g(x_1)-f(x_1), \ x_1\geq a\}$$ which is a $\beta$-strip at infinity, for some rational number $\beta\leq 1$, according to Proposition \ref{prop: triangle-f}.

\medskip
\noindent{\it Case $n>2$}.	In this case, $X$ can be either a graph or a band. First, let $X$ be a graph. Since $X$ is a graph of an outer Lipschitz function on a $\beta$-strip at infinity, we get $X$ itself is a $\beta$-strip at infinity ($\beta\leq 1$). Now, let us consider $X$ is a band 
$$X=\{(x',x_n)\in\R^{n-1}\times\R \ \colon \ f(x')<x_n<g(x'); \ x'\in X'\}$$
\noindent where $X'\subset\R^{n-1}$ is a (1-dimensional) standard L-regular cell in $\R^{n-1}$ with constant $C$, $f,g\colon X'\rightarrow\R$ are $C^1$ semialgebraic functions such that 
$$f(x')<g(x') \quad \forall \ x'\in X'$$ and $$ |df(x')|\leq C, \ |dg(x')|\leq C \quad \forall \ x'\in X'.$$ 

Since $X'$ is 1-dimensional, unbounded and has only one end, there exists a $C^1$-semialgebraic parametrization (outer lipeomorphism) $\gamma\colon(a,\infty)\rightarrow X'$; hence the closure $\overline{X}$ is semialgebraicly  outer lipeomorphic to the set $$\{(t,s)\in\R\times\R \ \colon \ f\circ\gamma(t)\leq s \leq g\circ\gamma(t), \ t\geq a\}$$ which implies $\overline{X}$ is a $\beta$-strip at infinity for some rational number $\beta\leq 1$.
\end{proof}

It is proved in \cite{KP}, more precisely, see Proposition 1.4 in \cite{KP} that any semialgebraic subset $X\subset\R^n$ can be stratified by L-regular cells in $\R^n$ with a constant $C=c(n)>0$. Now, we are ready to prove Theorem \ref{thm:beta_ends}.

\begin{proof}[Proof of Theorem \ref{thm:beta_ends}]
	Let $R>0$ be a sufficient large radius such that the connected components of $S\setminus B(0,R)$ are the ends of $S$. Let $X=S_i$ be one of those ends. As we mentioned above, we have a stratification $\displaystyle X=\bigcup_{i=1}^r C_i$ such that each stratum $C_i$ is a L-regular cell in $\R^n$ with constant $C=c(n)>0$. By taking $R>0$ large enough, one may suppose that all 2-dimensional strata of $X$ are unbounded (hence all of them have only one end). Then, let $C_{i_1},\dots,C_{i_k}$ be the 2-dimensional strata of $X$. It follows from Proposition \ref{prop:regular-cell} that the closure $\overline{C_{i_j}}$ of each cell $C_{i_j}$ is a $\beta_j$-strip at infinity for some rational number $\beta_j\leq 1$. Since $$X=\bigcup_{j=1}^k\overline{C_{i_j}}$$ and, by topological restrictions, the family $\overline{C_{i_1}},\dots,\overline{C_{i_k}}$ satisfies the assumptions of Proposition \ref{prop:glue-tube}, we get $X$ is a $\beta$-tube where $\beta=\max\{\beta_1,\dots,\beta_k\}$.
\end{proof}

\section{Classification of semialgebraic surfaces }\label{sec:classification}
In this Section, we are going to present a classification of all semialgebraic surfaces with isolated singularities.

\begin{remark}\label{rem:horn_exponent}
	Let $X\subset\R^n$ be a closed 2-dimensional semialgebraic set. According to the notion of topological regular points, we can read the Theorem \ref{thm:birbrair} (Theorem of Birbrair) in the following way: if  $p\in X$ is a topological regular point, then there exist a neighborhood $V\subset  X$ and a semialgebraic inner lipeomorphism $\phi\colon V\rightarrow H_{\beta}$; $\phi(p)=(0,0,0)$, where $\beta\geq 1$ is a rational number and $$H_{\beta}=\{(x,y,z)\in\R^3 \ \colon \ x^2+y^2=z^{2\beta} \ \mbox{and} \ z\geq 0\}.$$ The rational number $\beta$ is called the {\bf horn exponent} of $X$ at $p$. Notice that, it also follows from Theorem of Birbrair that a point $p\in X$ is inner Lipschitz regular if, and only if, the horn exponent of $X$ at $p$ is equal to $1$. 
\end{remark}

\begin{definition}\label{def:symbols_code}
Let $X\subset\R^n$ be a semialgebraic surface with isolated inner Lipschitz singularities. Let us consider the following symbols:
\begin{enumerate}
    \item [i)] For $p\in {\rm Sing}_{inLip}(X) $, $\ell(X,p)$ denotes the number of connected components of the link of $X$ at $p$;
	\item [ii)] We can consider a sufficient large radius $R>0$ and a small enough radius $\rho >0$ such that 
	$$
	X'=(X\cap \overline{B(0,R)})\setminus \bigg\{ B(x_1,\rho)\cup\cdots\cup B(x_s,\rho) \bigg\}
	$$
	is a topological surface with boundary and its topological type does not depend on $R$ and $\rho$. Thus, we define   
	$$\theta(X)= 
	\begin{cases}
	\ \ 1, \ \mbox{if} \ X' \ \mbox{is orientable} \\
	-1, \ \mbox{if} \ X' \ \mbox{is not orientable}.
	\end{cases}
	$$
	\item[iii)] $g(X)$ is the genus of $X'$;
	\item [iv)] For each $p\in X $, there is $r>0$ such that 
	$$X\cap B(p,r)=\bigcup\limits_{i=1}^{\ell(X,p)}X_i$$
	and each $X_i$ is a topological surface. Let $\beta_i$ be the horn exponent of $X_i$ at $p$ (see Remark \ref{rem:horn_exponent}). By reordering the indices, if necessary, we assume that $\beta_1\leq \beta_2\leq \cdots \leq \beta_{\ell(X,p)}$. In this way, we define $\beta(X,p)=(\beta_1, \beta_2, \cdots, \beta_{\ell(X,p)})$. 
	\item[v)] $e(X)$ is the number of ends of $X$, and if $E_1,\dots,E_{e(X)}$ are the ends of $X$, then denote by $\beta_i$, the tube exponent of $E_i$, the only rational number smaller than or equal to 1 such that $E_i$ is a $\beta_i$-tube. By reordering the indices, if necessary, we assume that $\beta_1\leq \beta_2\leq \cdots \leq \beta_{e(X)}$. In this way, we define $\beta(X,\infty)=(\beta_1, \beta_2,...,\beta_{e(X)})$. 
\end{enumerate}
\end{definition}

\begin{definition}[Inner Lipschitz code]\label{def:inner_code}
	Let $X\subset\R^n$ be a semialgebraic surface with isolated inner Lipschitz singularities. Let $S=\{p_1,...,p_k\}\subset X$ be a finite subset such that ${\rm Sing}_{inLip}(X)\subset S$ and let $\sigma\colon S\to \tilde S$ be a bijection for some subset $\tilde S$ in some Euclidean space.
	\begin{itemize}
	 \item If ${\rm Reg}_{inLip}(X)$ is a connected set, then the collection of symbols 
	$$\bigg\{ \theta(X),g(X),\beta(X,\infty),\{(\sigma(p);\beta(X,p))\}_{p\in S} 
	\bigg\}$$
	is called the {\bf inner Lipschitz code of $X$ w.r.t. $\sigma$} and we denote it by ${\rm Code}_{inLip}(X,\sigma)$. The collection of symbols 
	$$\bigg\{ \theta(X),g(X),\beta(X,\infty),\{\beta(X,p)\}_{p\in S} 
	\bigg\}$$
	is called the {\bf inner Lipschitz code of $X$} and we denote it by ${\rm Code}_{inLip}(X)$;
	\item For the general case, let $C_1,...,C_r$ be the closure of the connected components of ${\rm Reg}_{inLip}(X)$. The collection of inner Lipschitz codes 
	$$\bigg\{ {\rm Code}_{inLip}(C_1,\sigma|_{C_1\cap S}),\cdots, {\rm Code}_{inLip}(C_r,\sigma|_{C_r\cap S}) 
	\bigg\}$$
	is called the {\bf inner Lipschitz code of $X$ w.r.t. $\sigma$} and we also denote it by ${\rm Code}_{inLip}(X,\sigma)$. When $S={\rm Sing}_{inLip}(X)$ and $\sigma$ is the identity, we only denote ${\rm Code}_{inLip}(X,\sigma)$ by ${\rm Code}_{inLip}(X)$ and we also call it the {\bf inner Lipschitz code of $X$}.
	\end{itemize}
\end{definition}

\begin{example}
	Let us see the inner Lipschitz code of some well-known semialgebraic topological surfaces.
	\begin{enumerate}
		\item[a)] Right cylinder: $\{1,0,(0,0),\emptyset\}$;
		\item[b)] Unbounded Moebius band $\{(x,y,u,v)\in \R^4: \ x^2+y^2=1, \  (u^2-v^2)y=2uv x \}$: $\{-1,0,1,\emptyset\}$;
		\item[c)] Global $\beta$-horn in $\R^3$; $\beta\geq 1$: $\{1,0,1,\{\beta\}\}$;
		\item[d)] $\{(z,w)\in\mathbb{C}^2 \ \colon \ z^2=w(w-a)(w-b)\}$; $a,b\neq 0$ and $a\neq b$: $\{1,1,(1,1,1),\emptyset \}$;
		\item[e)] Paraboloid in $\R^3$: $\{1,0,1/2,\emptyset\}$
		\item[f)] Torus: $\{1,1,\emptyset,\emptyset\}$
		\item[g)] Klein bottle: $\{-1,1,\emptyset,\emptyset\}$
		\item[h)] Edge of two spheres $\{(x,y,z\in\R^3;((x-1)^2+y^2+z^2-1)((x+1)^2+y^2+z^2-1)=0\}$ : $\{\{1,0,\emptyset,\{((0,0,0);1)\}, \{1,0,\emptyset,\{((0,0,0);1)\}\}$
		\item[i)] Cayley surface $\{(x,y,z)\in\R^3;x^2+y^2+z^2-2xyz=1\}$ (see Figure 1): $\{\{1,0,1,(p_1;1)\}$, $\{1,0,1,(p_2;1)\}$, $\{1,0,1,(p_3;1)\}$, $\{1,0,1,(p_4;1)\}$, $\{1,0,\emptyset,$ $\{(p_1;1),(p_2;1),(p_3;1),(p_4;1)\}\}\}$.
	\end{enumerate}
\end{example}

\begin{figure}[H]
    \centering \includegraphics[scale=0.04]{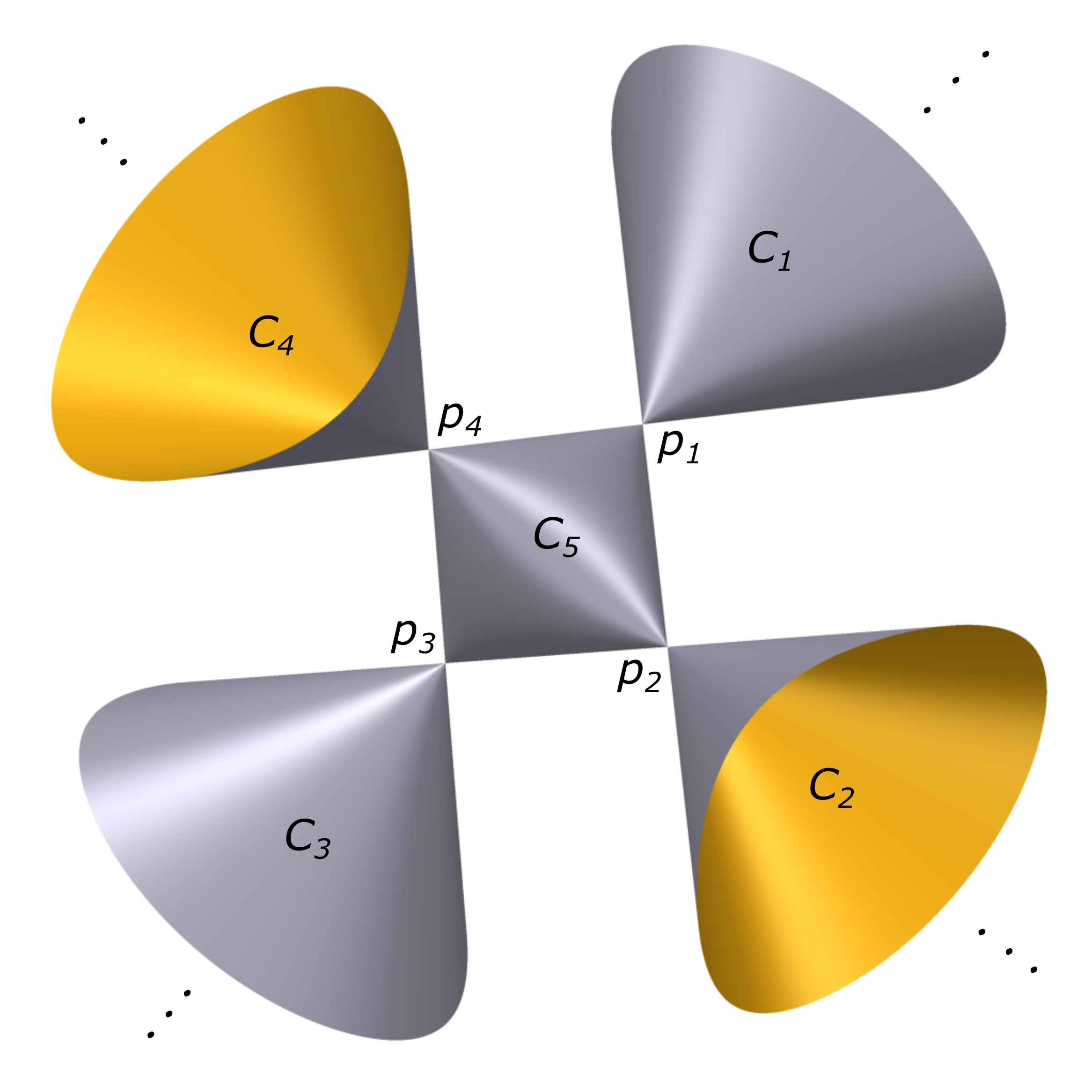}
    \caption{Decomposition of the Cayley surface $\{(x,y,z)\in\R^3;x^2+y^2+z^2-2xyz=1\}$.}\label{figure1}
\end{figure}

\begin{definition}
Let $X$ and $Y$ be two semialgebraic sets. We say that ${\rm Code}_{inLip}(X)$ and ${\rm Code}_{inLip}(Y)$ are equivalent if one of the following items holds true:
\begin{enumerate}
 \item ${\rm Reg}_{inLip}(X)$ and ${\rm Reg}_{inLip}(Y)$ are connected sets and ${\rm Code}_{inLip}(X)={\rm Code}_{inLip}(Y)$;
 \item ${\rm Reg}_{inLip}(X)$ and ${\rm Reg}_{inLip}(Y)$ are disconnected sets and ${\rm Code}_{inLip}(X,\sigma)={\rm Code}_{inLip}(Y)$ for some bijection $\sigma\colon {\rm Sing}_{inLip}(X)\to {\rm Sing}_{inLip}(Y)$.
\end{enumerate}

\end{definition}

\begin{theorem}\label{thm:class_surf}
	Let $X\subset\R^n$ and $Y\subset\R^m$ be  semialgebraic surfaces with isolated inner Lipschitz singularities. Then, $X$ and $Y$ are inner lipeomorphic if, and only if, their inner Lipschitz code are equivalent.
\end{theorem}

\begin{proof} Of course, we have the inner Lipschitz code is an inner Lipschitz invariant in the sense: if $X$ and $Y$ are inner lipeomorphic, then their codes are equivalent. From another hand, let us suppose that the inner Lipschitz codes of $X$ and $Y$ are equivalent. 

Let us assume, initially, that ${\rm Reg}_{inLip}(X)$ (and, consequently, ${\rm Reg}_{inLip}(Y)$) is a connected set.

Let us denote by $E_1^X,\dots,E_{e}^X$ the ends of $X$, with respective tube exponents $\beta_1(X)\leq \dots\leq\beta_{e}(X)$, and $E_1^Y,\dots,E_{e}^Y$ the ends of $Y$, with respective tube exponents $\beta_1(Y)\leq \dots\leq\beta_{e}(Y)$. Also, let us denote by $x_1,\dots,x_s$ the inner Lipschitz singularities of $X$, with respective horn exponents $\beta(X,x_1),\dots,\beta(X,x_s)$, and $y_1,\dots,y_s$ the inner Lipschitz singularities of $Y$, with respective horn exponents $\beta(Y,y_1),\dots,\beta(Y,y_s)$. So, we are assuming that $\theta(X)=\theta(Y)$, $g(X)=g(Y)$, $\beta(X,\infty)=\beta(Y,\infty)$ and $\beta(X,x_j)=\beta(Y,y_j)$, $j=1,\dots,s$. Then, we can consider $R>0$ a sufficient large radius and $\rho >0$ a small enough radius such that, for each $i\in \{1,\dots,e\}$ and $j\in \{1,\dots,s\}$, there exist semialgebraic inner lipeomorphisms  
$$ 
h_i \colon E_i^X\setminus B(0,R)\rightarrow E_i^Y\setminus B(0,R)\quad \mbox{and}\quad g_{j}\colon X\cap \overline{B(x_j,\rho)}\to Y\cap \overline{B(y_j,\rho)}.
$$

In fact, the existence of the $h_i$'s follows from Theorem \ref{thm:beta_ends}, and by writing $X\cap \overline{B(x_j,\rho)}=\bigcup\limits_{\ell=1}^{\ell(X,x_j)}X_{j\ell}$ (resp. $Y\cap \overline{B(y_j,\rho)}=\bigcup\limits_{\ell=1}^{\ell(Y,y_j)}Y_{j\ell}$) and $X_{j\ell}\cap X_{j\ell'}=\{x_j\}$ (resp. $Y_{j\ell}\cap Y_{j\ell'}=\{y_j\}$) whenever $\ell\not =\ell'$, by Theorem of Birbrair, there are inner lipeomorphims $g_{j\ell}\colon X_{j\ell}\to Y_{j\ell}$. So, we define $g_{j}\colon X\cap \overline{B(x_j,\rho)}\to Y\cap \overline{B(y_j,\rho)}$ by $g_{j}(z)=g_{j\ell}(z)$ whenever $z\in X_{j\ell}$.

Now, we consider the following Lipschitz surfaces with boundary
$$X'=(X\cap \overline{B(0,R)})\setminus \bigg\{ B(x_1,\rho)\cup\cdots\cup B(x_s,\rho) \bigg\}$$ and $$Y'=(Y\cap \overline{B(0,R)})\setminus \bigg\{ B(y_1,\rho)\cup\cdots\cup B(y_s,\rho) \bigg\},$$ and the following semialgebraic lipeomorphism $\kappa\colon\partial X'\rightarrow\partial Y'$ given by:
$$\kappa(z)=
\begin{cases} h_i(z), \ \mbox{if} \ z\in E_i^X; \ |z|=R \\
			g_{j}(z), \ \mbox{if} \ z\in X; \ |z-x_j|=\rho
\end{cases}.$$

Since $X'$ is orientable if, and only if, $Y'$ is orientable too, and $X'$ has the same genus and same number of boundary components as $Y'$, it follows from Proposition \ref{prop:lipeotopy}, maybe after changing the orientation of some $h_i$'s and $g_{j\ell}$'s, the following result.

\begin{lemma}\label{lemma:extension}
	 There exists a lipeomorphism $\Phi\colon X'\rightarrow Y'$ that extends $\kappa\colon\partial X'\rightarrow\partial Y'$
\end{lemma}

Finally, the mapping $F\colon X\rightarrow Y$ defined below is an inner lipeomorphism:

$$F(z)=
\begin{cases} h_i(z), \ \mbox{if} \ z\in E_i^X; \ |z|\geq R \\
g_j(z), \ \mbox{if} \ z\in X; \ |z-x_j|\leq \rho \\
\Phi(z), \mbox{if} \ z\in X'
\end{cases},$$  
which finishes the proof in this case.
 
Now, we have to consider the case that ${\rm Reg}_{inLip}(X)$ and ${\rm Reg}_{inLip}(Y)$ are disconnected sets. For this case, let $X_1,...,X_r$ (resp. $Y_1,...,Y_r$) be the closure of the connected components of ${\rm Reg}_{inLip}(X)$ (resp. ${\rm Reg}_{inLip}(X)$). We have assumed that ${\rm Code}_{inLip}(X)$ and ${\rm Code}_{inLip}(Y)$ are equivalent, then ${\rm Code}_{inLip}(X,\sigma)={\rm Code}_{inLip}(Y)$ for some bijection $\sigma\colon S={\rm Sing}_{inLip}(X)\to \tilde S={\rm Sing}_{inLip}(Y)$. By reordering the indices, if necessary, we may assume that ${\rm Code}_{inLip}(X_i,\sigma|_{X_i\cap S})={\rm Code}_{inLip}(Y_i, id_{\tilde S}|_{Y_i\cap \tilde S})$, $i=1,...,r$, where $id_{\tilde S}\colon \tilde S\to \tilde S$ is the identity mapping. 

For a closed semialgebraic set $A$ and $p\in A$, we have that $p\in {\rm Reg}_{inLip}(A)$ if and only if $\beta(A,p)=1$. Thus, fixed $i\in \{1,...,r\}$, for $S_i=(X_i\cap S)\setminus {\rm Sing}_{inLip}(X_i)$ and $\tilde S_i=(Y_i\cap \tilde S)\setminus {\rm Sing}_{inLip}(Y_i)$, we have $\sigma(S_i)=\tilde S_i$, and therefore ${\rm Code}_{inLip}(X_i)={\rm Code}_{inLip}(Y_i)$. By the first part of this proof, there is an inner lipeomorphism $F_i\colon X_i\to Y_i$. Moreover, we can take $F_i$ satisfying $F_i(p)=\sigma(p)$ for all $p\in X_i\cap S$. Thus, the mapping $F\colon X\rightarrow Y$, defined by $F(z)=
F_i(z)$ whenever $z\in X_i$, is an inner lipeomorphism, which finishes the proof.
\end{proof}

From now on, we start to list some consequences of Theorem \ref{thm:class_surf} and its proof.

\begin{figure}[H]
    \centering \includegraphics[scale=0.5]{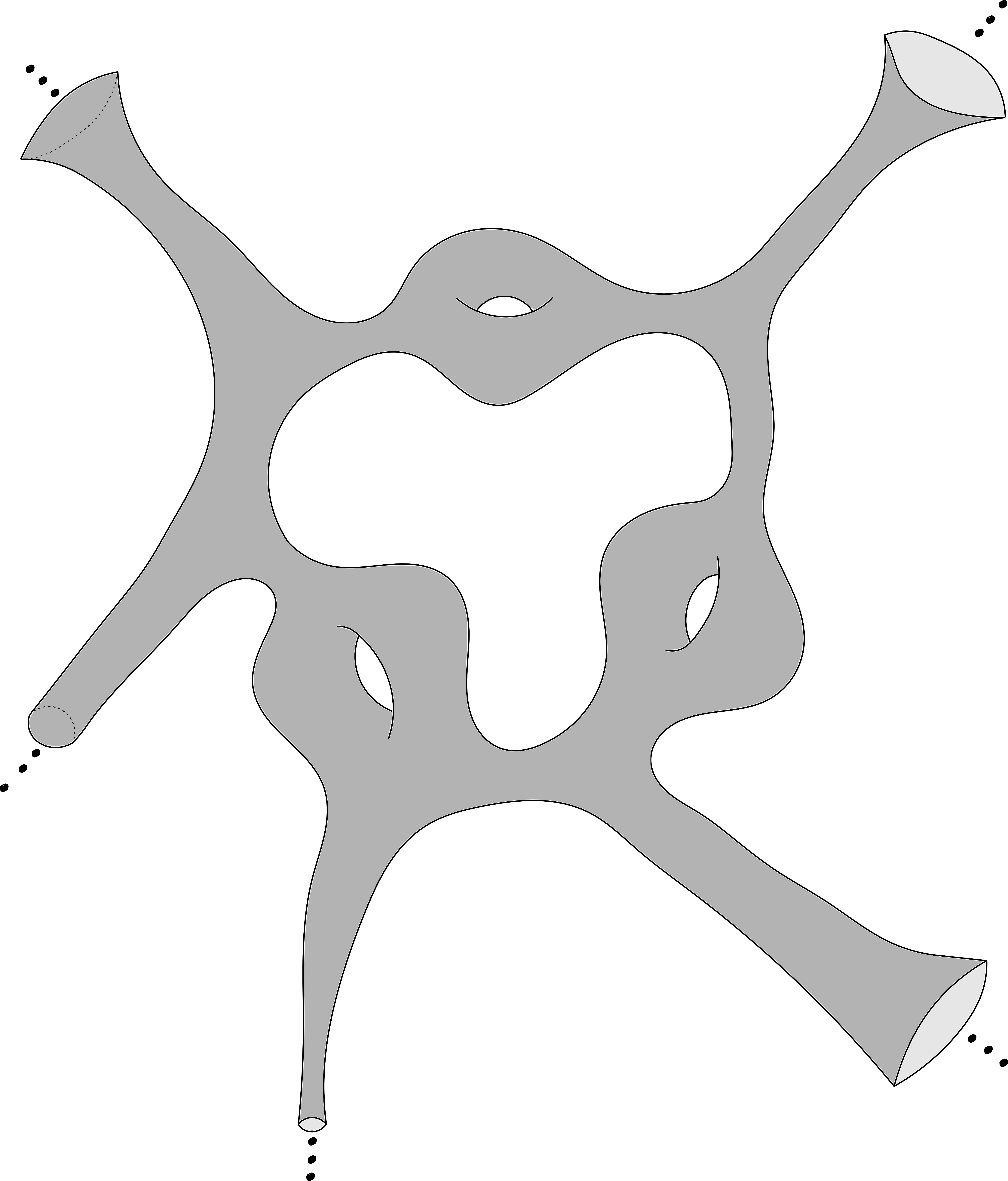}
    \caption{An oriented Nash surface with 5 ends and genus 4.}\label{figure2}
\end{figure}

The first consequence is a classification of the Nash surfaces, even for unbounded Nash surfaces as in Figure \ref{figure2}.

\begin{corollary}\label{cor:classf_nash_surfaces}
Let $N_1,N_2\subset \R^n$ be two Nash surfaces. Then, the following statements are equivalent:
\begin{itemize}
 \item [(1)] $N_1$ and $N_2$ are homeomorphic and $\beta(N_1,\infty)=\beta(N_2,\infty)$;
 \item [(2)] $N_1$ and $N_2$ are inner lipeomorphic;
 \item [(3)] $\theta(N_1)=\theta(N_2)$, $g(N_1)=g(N_2)$ and $\beta(N_1,\infty)=\beta(N_2,\infty)$.
\end{itemize} 
\end{corollary}

\begin{remark}
Since properly embedded smooth surfaces in $\R^3$ are orientable, we obtain that two Nash surfaces $N_1$ and $N_2$ in $\R^3$ are inner lipeomorphic if and only if $g(N_1)=g(N_2)$ and $\beta(N_1,\infty)=\beta(N_2,\infty)$.
\end{remark}

In fact, we obtain a stronger result than Corollary \ref{cor:classf_nash_surfaces}, since we can present normal forms for the classification presented in \ref{cor:classf_nash_surfaces}. In order to that, for $\theta \in \{-1,1\}$ and $g\in \mathbb{N}$, let $N(\theta, g)\subset \R^5$ be a compact Nash surface such that $\theta(N(\theta, g))=\theta$ and $g(N(\theta, g))=g$. For a positive integer number $e$ and $\beta=(\beta_1,...,\beta_e)\in \mathbb{Q}$ such that $\beta_1\leq \beta_2\leq ....\leq \beta_e\leq 1$, we remove $e$ distinct points of $N(\theta, g)$, let us say $x_1,...,x_e\in N(\theta, g)$, and we define $F\colon N(\theta, g)\setminus \{x_1,...,x_e\}\to \R^{6e}$ given by 
$$
F(x)=(\frac{x-x_1}{\|x-x_1\|^{1+\beta_1}},\|x-x_1\|^{-1},\frac{x-x_2}{\|x-x_2\|^{1+\beta_2}},\|x-x_2\|^{-1}, ..., \frac{x-x_e}{\|x-x_e\|^{1+\beta_e}},\|x-x_e\|^{-1}).
$$ 
We denote the image of $F$, which is a Nash surface, by $N(\theta, g,\beta)$. We also define $N(\theta, g,\emptyset)=N(\theta, g)$. Note that $\theta(N(\theta, g,\beta))=\theta$, $g(N(\theta, g,\beta))=g$ and $\beta(N(\theta, g,\beta),\infty)=\beta$. Thus, $N(\theta, g,\beta)$ is well defined up to inner lipeomorphisms, and we obtain the following:

\begin{corollary}\label{cor:normal_forms_classf_nash_surfaces}
Let $N\subset \R^n$ be a Nash surface. Then, $N(\theta(N),g(N),\beta(N,\infty))$ and $N$ are inner lipeomorphic. 
\end{corollary}

\begin{corollary}
Let $M_1,M_2\subset \R^3$ be two connected properly embedded minimal surfaces with finite total curvature. Then, the following statements are equivalent:
\begin{itemize}
 \item [(1)] $M_1$ and $M_2$ are homeomorphic;
 \item [(2)] $M_1$ and $M_2$ are inner lipeomorphic;
 \item [(3)] $g(M_1)=g(M_2)$ and $e(M_1)=e(M_2)$.
\end{itemize}
\end{corollary}
\begin{proof}
Obviously, (2) implies (1), and (1) implies (3).

Thus, we only have to show that (3) implies (2). Let us assume $g(M_1)=g(M_2)$ and $e(M_1)=e(M_2)$.

Since the tangent cone at infinity of each end of a properly embedded minimal surfaces with finite total curvature is a plane, it follows from, for example, Lemma 1 in \cite{BelenkiiB:2005} that such an end is inner lipeomorphic to $\R^2$. By proof of Theorem \ref{thm:class_surf}, $M_1$ and $M_2$ are inner lipeomorphic.
\end{proof}

\begin{remark}\label{non-degen-tubes}
Let $X\subset \R^n$ be a closed semialgebraic surface which is a $\beta$-tube. Then, $\beta=1$ if and only if $\dim C(X,\infty)=2$.
\end{remark}

\begin{corollary}
Let $C_1,C_2\subset \mathbb{C}^2$ be two complex algebraic curves. Then, the following statements are equivalent:
\begin{itemize}
 \item [(1)] $C_1$ and $C_2$ are homeomorphic;
 \item [(2)] $C_1$ and $C_2$ are inner lipeomorphic;
 \item [(3)] If $X_1,...,X_r$ and $Y_1,...,Y_s$ are the irreducible components of $C_1$ and $C_2$, respectively, then there exist bijections  $\pi\colon \{1,...r\}\to \{1,...,s\}$ and $\sigma\colon {\rm Sing}_{inLip}(C_1)\to {\rm Sing}_{inLip}(C_2)$ such that $g(X_i)=g(Y_{\pi(i)})$, $e(X_i)=e(Y_{\pi(i)})$ and $\ell(X_i,p)=\ell(Y_{\pi(i)},\sigma(p))$ for all $p\in {\rm Sing}_{inLip}(C_1)$, $i=1,...,r$.
\end{itemize}
\end{corollary}
\begin{proof}
Obviously, (2) implies (1). 

We are going to show that (1) implies (2). Assume that there is a homeomorphism $h\colon C_1\to C_2$. Therefore, $g(C_1)=g(C_2)$ and $e(C_1)=e(C_2)$. Since the tangent cone at infinity of each end of complex algebraic curve is a complex line, by Remark \ref{non-degen-tubes}, we obtain that each such end is a $1$-tube. Thus $\beta(C_1,\infty)=\beta(C_2,\infty)$. 

It follows from the Birbrair Theorem that two irreducible germs of complex analytic are inner lipeomorphic. In particular, $p\in {\rm Sing}_{inLip}(C_1)$ if and only if $\ell(C_1,p)>1$. Thus, $h({\rm Sing}_{inLip}(C_1))={\rm Sing}_{inLip}(C_2)$ and $\sigma=h|_{{\rm Sing}_{inLip}(C_1)}$ is a bijection. Moreover, each irreducible germ of a complex analytic curve is the germ of a $1$-horn. Therefore, $\beta(C_1,p)=\beta(C_2,\sigma(p))$ for all $p\in {\rm Sing}_{inLip}(C_1)$. Then ${\rm Code}_{inLip}(C_1,\sigma)={\rm Code}_{inLip}(C_2)$. By Theorem \ref{thm:class_surf}, $C_1$ and $C_2$ are inner lipeomorphic.

In order to finish the proof, due to the comments made in this proof, we note that the item (3) is equivalent to say that ${\rm Code}_{inLip}(C_1,\sigma)$ and ${\rm Code}_{inLip}(C_2)$ are equivalent, which finishes the proof.
\end{proof}

%
%
%
%
%

\end{document}